\documentclass[english,11pt]{article}
\usepackage[latin1]{inputenc}
\usepackage{amsmath,amsthm,amssymb,babel}

\def\eq#1{(\ref{#1})}

\def\neweq#1{\begin{equation}\label{#1}}
\def\endeq{\end{equation}}

\def\ep{\varepsilon}
\def\la{\lambda}

\def\phi{\varphi}
\def\RR{{\mathbb R} }

\def\di{\displaystyle}
\def\ri{\rightarrow}

\newtheorem{theorem}{Theorem}[section]
\newtheorem{lem}{Lemma}[section]

\title{\sc On a class of sublinear
singular elliptic problems with convection term}
\author{Marius GHERGU and
Vicen\c tiu R\u ADULESCU\\
\small Department of Mathematics, University of Craiova, 200585
Craiova, Romania}
\date{}
\renewcommand{\theequation}{\arabic{section}.\arabic{equation}}
\catcode`@=11 \@addtoreset{equation}{section} \catcode`@=12
\begin{document}
\baselineskip16pt
\maketitle
\renewcommand{\theequation}{\arabic{section}.\arabic{equation}}
\catcode`@=11 \@addtoreset{equation}{section} \catcode`@=12

\bigskip
\noindent
{\large\bf Correspondence address:}\\ Vicen\c
tiu R\u adulescu\\ Department of Mathematics\\ University of
Craiova\\ 200585 Craiova, Romania\\ fax: +40-251.411688\\ E-mail:
{\tt vicentiu.radulescu@ucv.ro}

\newpage
\begin{abstract}
We establish several results related to existence, nonexistence or
bifurcation of positive solutions for the boundary value problem
$-\Delta u+K(x)g(u)+|\nabla u|^a=\la f(x,u)$ in 
$\Omega$, $u=0$ on $\partial
\Omega,$
where $\Omega\subset\RR^N$ $(N\geq 2)$ is a smooth
bounded domain,   $0<a\leq 2,$ $\la$ is a positive parameter, 
and $f$ is smooth and has a sublinear growth.
The main feature of this paper consists in the presence of the singular
nonlinearity $g$ combined with the convection term $|\nabla u|^a.$ 
Our approach takes into account both 
the sign of the potential $K$ and the decay rate around the 
origin of the singular nonlinearity $g$.
The proofs are based on  various techniques related to
the maximum principle for elliptic equations.
\end{abstract}

\section{Introduction and the main results}
Stationary problems involving singular nonlinearities,
 as well as the associated evolution equations, 
describe naturally several physical phenomena. 
At our best knowledge, the first study in this direction is due to Fulks and Maybee
\cite{ful}, who proved existence and uniqueness results by using a fixed point argument; moreover, 
they showed that solutions of the parabolic problem tend to the unique solution of the 
corresponding elliptic equation. A different approach (see \cite{coq,crt,stu}) consists in 
approximating the singular equation with a regular problem, where the standard techniques
(e.g., monotonicity methods) can be applied and then passing to the limit to obtain the solution 
of the original equation. 
Nonlinear singular boundary value problems  
arise in the context of chemical heterogeneous catalysts and chemical catalyst kinetics, in the 
theory of heat conduction in electrically conducting materials, singular minimal surfaces, as well 
as in the study of non-Newtonian fluids, boundary layer phenomena for viscous fluids (we refer for 
more details to \cite{bn,caf,cn1,cn,diaz,dmo}
and the more recent papers \cite{hai,her,her1,mea,sy1,sy2,z1}). We also point out that, due to the 
meaning of the unknowns (concentrations, populations, etc.), only the positive solutions are  
relevant in most cases. 

Let $\Omega$ be a smooth bounded domain in
$\RR^N$ ($N\geq 2$).
We are concerned in this paper with the
following boundary value problem
$$\left\{\begin{tabular}{ll}
$-\Delta u+K(x)g(u)+|\nabla u|^a=\la f(x,u)$ \quad & ${\rm
in}\ \Omega,$\\
$u>0$ \quad & ${\rm in}\ \Omega,$\\
$u=0$ \quad & ${\rm on}\ \partial\Omega,$\\
\end{tabular} \right. \eqno(1)_\la$$
where $\la>0,$ $0<a\leq 2$ and $K\in C^{0,\gamma}(\overline\Omega)$, $0<\gamma<1.$ 
Here $f:\overline{\Omega}\times[0,\infty)\rightarrow[0,\infty)$
is a H\"{o}lder continuous function
which is positive on
$\overline{\Omega}\times(0,\infty).$
We  assume that $f$ is nondecreasing with
respect to the second variable and is sublinear, that is,
\medskip

\noindent $\di(f1)\qquad$ the mapping $\di
(0,\infty)\ni
s\longmapsto\frac{f(x,s)}{s}\quad\mbox{is
nonincreasing for all}\;\, x\in\overline{\Omega};$
\medskip

\noindent $\di(f2)\qquad \lim_{s\ri 0^+}\frac{f(x,s)}{s}=+\infty\quad\mbox{and}\;\;
\lim_{s\rightarrow\infty}\frac{f(x,s)}{s}=0,\;\;\mbox{uniformly
for}\;\,x\in\overline{\Omega}.$
\medskip

\noindent We assume that $g\in
C^{0,\gamma}(0,\infty)$ is a nonnegative and
nonincreasing function satisfying
\medskip

\noindent $\di(g1)\qquad \lim_{s\ri 0^+}g(s)=+\infty.$
\medskip

Problem $(1)_\la$ has been considered in \cite{gr1}
in the absence of the gradient term $|\nabla u|^a$ and assuming that the singular 
term $g(t)$ behaves like $t^{-\alpha}$ around the origin, with $t\in (0,1)$. 
In this case it has been shown that the  
sign of the extremal values of $K$ plays a crucial role. 
In this sense, we have
proved in \cite{gr1} that if $K<0$ in $\overline\Omega,$ then problem $(1)_\la$ (with $a=0$) has a
unique solution in the class ${\cal E}=\{u\in C^2(\Omega)\cap
C({\overline{\Omega}});\,\,g(u)\in L^1(\Omega)\},$ for all
$\la>0.$ On the other hand, 
if $K>0$ in $\overline\Omega$, then there exists $\la^*$ such that
problem $(1)_\la$ has solutions in $\cal E$ if $\lambda>\la^*$ and no solution exists if 
$\lambda<\la^*$. 
The case where $f$ is asymptotically linear, $K\leq 0$, and $a=0$ has been discussed in 
\cite{crg}. 
In this case, a major role is played by 
$\lim_{s\ri\infty}f(s)/s=m>0.$
More precisely, there exists a solution (which is unique) $u_\la\in C^2(\Omega)\cap 
C^{1}(\overline\Omega)$
 if and only if $\la<\la^*:=\la_1/m.$ An additional result asserts that the mapping 
$(0,\la^*)\longmapsto
u_\lambda$  is increasing and  
$\lim_{\la\nearrow \la^*} u_\la=+\infty$ uniformly on 
compact subsets of $\Omega.$

Due to the singular character of our problem $(1)_\la,$ 
we cannot expect to have solutions in $C^2(\overline\Omega).$ 
We are seeking in this paper
classical solutions of $(1)_\la,$ that is, solutions $u\in
C^2(\Omega)\cap C(\overline\Omega)$ that verify $(1)_\la.$ 
Closely related to our problem is the following one, which has been considered 
in \cite{gr3}:
\neweq{edinburgh}
\left\{\begin{tabular}{ll}
$-\Delta u=g(u)+|\nabla u|^a+\la f(x,u)$ \quad & $\mbox{\rm in}\ \Omega,$\\
$u>0$ \quad & $\mbox{\rm in}\ \Omega,$\\
$u=0$ \quad & $\mbox{\rm on}\ \partial\Omega,$\\
\end{tabular} \right.
\endeq
where $f$ and $g$ verifies the above assumptions $(f1),$ $(f2)$ 
and $(g1).$ We have proved in \cite{gr3} that if 
$0<a<1$ then problem
\eq{edinburgh} has at least one classical solution for all
$\la\geq 0.$ In turn, if $1<a\leq 2,$ 
then \eq{edinburgh} has no solutions for large values of $\la>0.$

The existence results for our problem $(1)_\la$ are quite different 
to those of \eq{edinburgh} presented in \cite{gr3}. More exactly we 
prove in the present paper that problem 
$(1)_\la$ has at least one solution only when $\la>0$ is 
large enough and $g$ satisfies a naturally 
growth condition around the origin. 
We extend the results in \cite[Theorem 1]{bdd}, corresponding to $K\equiv 0,$ 
$f\equiv f(x)$ and $a\in[0,1).$

The main difficulty in the treatment of $(1)_\la$ is the 
lack of the usual maximal principle between super and
sub-solutions, due to the singular
character of the equation. To overcome it, we state an 
improved comparison principle that fit to our problem $(1)_\la$ (see Lemma \ref{l} below).   

Throughout this paper we assume that $f$ satisfies
assumptions $(f1)-(f2)$  and $g$ verifies condition $(g1).$

In our first result we assume that $K<0$ in $\Omega.$
Note that $K$ may vanish on $\partial\Omega$ which 
leads us to a competition on the boundary between the 
potential $K(x)$ and the singular term $g(u).$ 
We prove the following result.

\begin{theorem}\label{th0}
Assume that $K<0$ in $\Omega.$ 
Then, for all $\la>0$, problem $(1)_\la$ has at least one classical solution.
\end{theorem}

Next, we assume that $K>0$ in $\overline\Omega.$ In this case, 
the existence of a solution to
$(1)_\la$ is closely related to the decay rate around its 
singularity. In this sense, we prove that problem $(1)_\la$
has no solution, provided that $g$ has a ``strong" singularity at the origin. More precisely, we 
have

\begin{theorem}\label{th1}
Assume that   $K>0$ in $\overline\Omega$ and 
$\int^1_0 g(s)ds=+\infty$. Then
problem $(1)_\la$ has no classical solutions.
\end{theorem}

In the following result, assuming that $\int^1_0 g(s)ds<+\infty$, we 
show that problem $(1)_\la$ 
has at least one solution, provided that $\la>0$ is large enough. 
Obviously, the hypothesis $\int^1_0 g(s)ds<+\infty$ implies
the following Keller-Osserman type condition around the origin

\medskip
\noindent $\di(g3)\qquad \int^1_{0}\left(
\int^t_0g(s)ds\right)^{-1/2}dt<\infty.$

\medskip 

As proved by B\'enilan, Brezis and Crandall
\cite{bbc}, condition
$(g3)$ is equivalent to the {\it property of compact
support}, that is,
for every $h\in L^1(\RR^N)$ with compact support,
there exists
a unique $u\in W^{1,1}(\RR^N)$ with compact support
such that
$\Delta u\in L^1(\RR^N)$ and
$$-\Delta u+g(u)=h\qquad\mbox{a.e. in}\ \RR^N.$$

\begin{theorem}\label{th2}
Assume that $K>0$ in $\overline\Omega$ and $\int^1_0 g(s)ds<+\infty.$
Then there exists $\la^*>0$ such that problem
$(1)_\la$ has at least one classical solution if
$\la>\la^*$ and no solution exists if $\la<\la^*.$
\end{theorem}

In the next section we establish a general comparison result between sub and super-solutions. 
Sections 3, 4 and 5 are devoted to the proofs of the above theorems.

\section{A comparison principle}
A very useful auxiliary result is the following comparison principle 
that improves Lemma 3 in \cite{sy1}. The proof uses some ideas from 
Shi and Yao \cite{sy1}, that goes back to the pioneering work by Brezis and Kamin \cite{bk}.

\begin{lem}\label{l}
Let
$\Psi:\overline{\Omega}\times(0,\infty)\rightarrow\RR$
be a continuous function
such that the mapping  $\di(0,\infty)\ni
s\longmapsto\frac{\Psi(x,s)}{s}$ is strictly decreasing
at
each $x\in\Omega.$ Assume that there exists  $v$,
$w\in C^2(\Omega)\cap C({\overline{\Omega}})$ such that

$(a)\qquad \Delta w+\Psi(x,w)\leq 0\leq \Delta v+\Psi(x,v)$
in $\Omega;$

$(b)\qquad v,w>0$ in $\Omega$ and $v\leq w$
on $\partial\Omega;$

$(c)\qquad\Delta v\in L^1(\Omega)\;\mbox{ or }\;
\Delta w\in L^1(\Omega).$\\
Then $v\leq w$ in $\Omega.$
\end{lem}

\noindent {\bf Proof.} We argue by contradiction and assume that
$v\geq w$ is not true in $\Omega.$ 
Then, we can find $\ep_0,\delta_0>0$ and a 
ball $B\subset\subset\Omega$ such that
$v-w\geq \ep_0$  in $B$ and
\neweq{shi2}
\di \int_Bvw\left(\frac{\Psi(x,w)}{w}-\frac{\Psi(x,v)}{v}\right)
dx\geq \delta_0.
\endeq
The case $\Delta v\in L^1(\Omega)$ was presented in 
\cite[Lemma 3]{sy1}. 
Let us assume now that $\Delta w\in L^1(\Omega)$ and set
$M=\max\{1,\|\Delta w\|_{L^1(\Omega)}\}$, 
$\ep=\min\left\{1,\ep_0,2^{-2}\delta_0/M\right\}.$
Consider  a nondecreasing function $\theta\in C^1(\RR)$ 
such that $\theta(t)=0,$ if  
$t\leq 1/2,$ $\theta(t)=1,$ if $t\geq 1,$ and 
$\theta(t)\in(0,1)$ if $t\in(1/2,1).$ Define
$$\di \theta_\ep(t)=\theta\left(\frac{t}{\ep}\right),\quad t\in\RR.$$
Since $w\geq v$ on $\partial\Omega,$ we can find a smooth 
subdomain $\Omega^*\subset\subset\Omega$ such that
$$B\subset\Omega^*\quad\mbox{ and }\;
v-w<\frac{\ep}{2}\;\mbox{ in }\,\Omega\setminus\Omega^*.$$ 
Using the hypotheses (a) and (b) we deduce
\neweq{shi3}
\di \int_{\Omega^*}(w\Delta v-v\Delta w)\theta_\ep(v-w)dx
\geq\int_{\Omega^*}
vw\left(\frac{\Psi(x,w)}{w}-\frac{\Psi(x,v)}{v}\right)
\theta_\ep(v-w)dx.
\endeq
By \eq{shi2} we have
$$\begin{tabular}{ll}
$\di\int_{\Omega^*}
vw\left(\frac{\Psi(x,w)}{w}-\frac{\Psi(x,v)}{v}\right)
\theta_\ep(v-w)dx$&$\di\geq 
\int_{B}
vw\left(\frac{\Psi(x,w)}{w}-\frac{\Psi(x,v)}{v}\right)
\theta_\ep(v-w)dx$\\
&$\di=\int_{B}
vw\left(\frac{\Psi(x,w)}{w}-\frac{\Psi(x,v)}{v}\right)dx$\\
&$\geq \delta_0.$\\
\end{tabular}$$
To raise a contradiction we need only to prove that the left-hand 
side in \eq{shi3} is smaller than 
$\delta_0.$ For this purpose, we define
$$\di \Theta_\ep(t)=\int_0^ts\theta'_\ep(s)ds,\quad t\in\RR.$$
It is easy to see that
\neweq{shi4}
\Theta_\ep(t)=0, \;\mbox{ if }\,t<\frac{\ep}{2}\quad\mbox{ and }\,\;
0\leq \Theta_\ep(t)\leq 2\ep,\;\mbox{ for all }t\in\RR.
\endeq
Now, using the Green theorem, we evaluate the left-hand side of \eq{shi3}:
$$\begin{tabular}{ll}
&$\di \int_{\Omega^*}(w\Delta v-v\Delta w)\theta_\ep(v-w)dx$\\
$=$&$\di 
\int_{\partial\Omega^*}w\theta_\ep(v-w)\frac{\partial v}
{\partial n}d\sigma-
\int_{\Omega^*}(\nabla w\cdot\nabla v)\theta_\ep(v-w)dx$\\
&$\di-\int_{\Omega^*}w\theta'_\ep(v-w)\nabla v
\cdot\nabla( v-w)dx-
\int_{\partial\Omega^*}v\theta_\ep(v-w)\frac{\partial w}
{\partial n}d\sigma$\\
&$\di+\int_{\Omega^*}(\nabla w\cdot \nabla v)\theta_\ep(v-w)dx+
\int_{\Omega^*}v\theta'_\ep(v-w)\nabla w\cdot\nabla(v-w)dx$\\
$=$&$\di\int_{\Omega^*}\theta'_\ep(v-w)(v\nabla w-w\nabla v)
\cdot\nabla( v-w)dx.$\\
\end{tabular}$$
The above relation can also be rewritten as
$$\begin{tabular}{ll}
$\di \int_{\Omega^*}(w\Delta v-v\Delta w)\theta_\ep(v-w)dx=$&
$\di\int_{\Omega^*}w\theta'_\ep(v-w)\nabla(w-v)
\cdot\nabla(v-w)dx$\\
&$\di+\int_{\Omega^*}(v-w)\theta'_\ep(v-w)\nabla w
\cdot\nabla(v-w)dx.$\\
\end{tabular}$$
Since $\di 
\int_{\Omega^*}w\theta'_\ep(v-w)\nabla(w-v)
\cdot\nabla(v-w)dx\leq 0,$ the last equality yields
$$\di \int_{\Omega^*}(w\Delta v-v\Delta w)\theta_\ep(v-w)dx\leq 
\di\int_{\Omega^*}(v-w)\theta'_\ep(v-w)\nabla w
\cdot\nabla(v-w)dx,$$ 
that is,
$$\di \int_{\Omega^*}(w\Delta v-v\Delta w)\theta_\ep(v-w)dx\leq 
\di\int_{\Omega^*}\nabla w\cdot \nabla(\Theta_\ep(v-w))dx.$$ 
Again by Green's first formula and by \eq{shi4} we have
$$\begin{tabular}{ll}
$\di \int_{\Omega^*}(w\Delta v-v\Delta w)\theta_\ep(v-w)dx$&$\di\leq
\di\int_{\partial\Omega^*}\Theta_\ep(v-w)
\frac{\partial v}{\partial n}d\sigma-
\di\int_{\Omega^*}\Theta_\ep(v-w)\Delta wdx$\\
&$\di\leq  -
\int_{\Omega^*}\Theta_\ep(v-w)\Delta wdx\leq 
2\ep\int_{\Omega^*}|\Delta w|dx$\\
&$\di\leq 2\ep M<\frac{\delta_0}{2}.$\\
\end{tabular}$$
Thus, we have obtained a contradiction. 
Hence $v\leq w$ in $\Omega$ and the proof of Lemma \ref{l} is now complete.
\qed

\section{Proof of Theorem \ref{th0}}
We need the following auxiliary result, which is proved in \cite{sy2}.

\begin{lem}\label{l2}
Let
$\Psi:\overline{\Omega}\times(0,\infty)\rightarrow\RR$
be a H\"{o}lder continuous function which satisfies
\smallskip

\noindent $\di(A1)\qquad \limsup _{s\rightarrow
+\infty}
\left(s^{-1}\max_{x\in\overline{\Omega}}\Psi(x,s)\right)<\la_1;$
\smallskip

\noindent $(A2)\qquad$ for each $t>0,$ there exists
a constant $D(t)>0$ such that
$$\Psi(x,r)-\Psi(x,s)\geq -D(t)(r-s),\quad\mbox{for}\;\;x\in
{\overline{\Omega}}\,
\;\;\mbox{and}\;\;\;r\geq s\geq t;$$
\smallskip

\noindent $(A3)\qquad$ there exist $\eta_0>0$
and an open subset
$\Omega_0\subset\Omega$ such that
$$\di\min_{x\in\overline{\Omega}}\Psi(x,s)\geq
0\quad\mbox{for}\;x\in(0,\eta_0),$$
and
$$\di \lim_{s\downarrow
0}\frac{\Psi(x,s)}{s}=+\infty\quad\mbox{uniformly for}
\;\,x\in\Omega_0.$$

Then the problem
\begin{equation}\label{doidoi}
 \left\{\begin{tabular}{ll}
$-\Delta u=\Psi(x,u)$ \quad & ${\rm in}\ \Omega,$\\
$u> 0$ \quad & ${\rm in}\ \Omega,$\\
$u=0$\quad & ${\rm on}\ \partial\Omega,$\\
\end{tabular} \right.
\end{equation}
has at least one classical solution $u\in
C^{2}(\Omega)\cap C(\overline{\Omega}).$
\end{lem}
Fix $\la>0.$ Obviously, $\Psi(x,s)=\la f(x,s)-K(x)g(s)$ satisfies the hypotheses 
in Lemma \ref{l2}
since $K<0$ in $\Omega.$ Hence, there exists a solution $\overline u_\la$ of the 
problem 
$$
\left\{\begin{tabular}{ll}
$-\Delta u=\la f(x,u)-K(x)g(u)$ \quad & ${\rm in}\ \Omega,$\\
$u>0$ \quad & ${\rm in}\ \Omega,$\\
$u=0$ \quad & ${\rm on}\ \partial\Omega.$\\
\end{tabular}\right. 
$$
We observe that $\overline u_\la$ is a super-solution of  problem $(1)_\la.$
To find a sub-solution, let us denote 
$$\di p(x)=\min\{\la f(x,1);-K(x)g(1)\}, \quad x\in\overline\Omega.$$
Using the monotonicity of $f$ and $g,$ we observe that
$p(x)\leq \la f(x,s)-K(x)g(s)$ for all $(x,s)\in\Omega\times (0,\infty).$ 
We now consider the problem
\neweq{subsol}
\left\{\begin{tabular}{ll}
$-\Delta v+|\nabla v|^a=p(x)$ \quad & ${\rm in}\ \Omega,$\\
$v=0$ \quad & ${\rm on}\ \partial\Omega.$\\
\end{tabular}\right. \endeq
First, we observe that $v=0$ is a sub-solution of \eq{subsol} while
$w$ defined by
$$\left\{\begin{tabular}{ll}
$-\Delta w=p(x)$ \quad & ${\rm in}\ \Omega,$\\
$w=0$ \quad & ${\rm on}\ \partial\Omega,$\\
\end{tabular}\right. $$
is a super-solution. Since $p>0$ in $\Omega$ we deduce that $w\geq 0$ in $\Omega.$
Thus, the problem \eq{subsol} has at least one 
classical solution $v.$ We claim that 
$v$ is positive in $\Omega.$ Indeed, if $v$ has a 
minimum in $\Omega,$ say at $x_0,$ then
$\nabla v(x_0)=0$ and $\Delta v(x_0)\geq 0.$ Therefore
$$\di 0\geq -\Delta v(x_0)+|\nabla v|^a(x_0)=p(x_0)>0,$$
which is a contradiction. Hence $\min_{x\in\overline\Omega} v=
\min_{x\in\partial\Omega} v=0,$ that is, $v>0$ in $\Omega.$ 
Now $\underline u_\la=v$ is a 
sub-solution of $(1)_\la$ and we have
$$\di -\Delta \underline u_\la=p(x)\leq \la f(x,\overline u_\la)-K(x)g(\overline u_\la)=
-\Delta \overline u_\la\quad\mbox{ in }\,\Omega.$$
Since $\underline u_\la=\overline u_\la=0$ on $\partial\Omega,$ from the above relation
we may conclude that $\underline u_\la\leq \overline u_\la$ in $\Omega$ and so, 
there exists at
least one classical solution for $(1)_\la.$ 
The proof of Theorem \ref{th0} is now complete.
\qed

\section{Proof of Theorem \ref{th1}}
We  give a direct proof, without using any change of variable, as in \cite{z1}. 
Let us assume that there exists $\la>0$ such that the problem $(1)_\la$ has a classical 
solution $u_\la.$ Since
$f$ satisfies $(f1)$ and $(f2)$,  
we deduce by Lemma \ref{l2} that for all $\la>0$
there exists $U_\la\in C^2(\overline\Omega)$  such that
\begin{equation}\label {UU}
 \left\{\begin{tabular}{ll}
$-\Delta U_\la=\la f(x,U_\la)$ \quad & ${\rm in}\
\Omega,$\\
$U_\la>0$ \quad & ${\rm in}\ \Omega,$\\
$U_\la=0$ \quad & ${\rm on}\ \partial\Omega.$\\
\end{tabular} \right.
\end{equation}
Moreover, there exist $c_1,c_2>0$ such that
\neweq{udst}
c_1\,\mbox{dist}\,(x,\partial\Omega)\leq U_{\la}(x)\leq
c_2\, \mbox{dist}\,(x,\partial\Omega)\quad\mbox{for all}
\;x\in\Omega.
\endeq
Consider the perturbed problem
\begin{equation}\label {epsilon}
 \left\{\begin{tabular}{ll}
$-\Delta u+K_*g(u+\ep)=\la f(x,u)$ \quad &
${\rm in}\ \Omega,$\\
$u>0$ \quad & ${\rm in}\ \Omega,$\\
$u=0$ \quad & ${\rm on}\ \partial\Omega,$\\
\end{tabular} \right.
\end{equation}
where $K_*=\min_{x\in\overline\Omega}K(x)>0.$
It is clear that $u_{\la}$ and
$U_{\la}$ are respectively sub and super-solution of (\ref{epsilon}).
Furthermore, we have
$$\di \Delta U_\la+f(x,U_\la)\leq 0\leq \Delta u_\la+f(x,u_\la)\quad\mbox{ in }\Omega,$$
$$\di U_\la,u_\la>0\quad\mbox{ in }\Omega,$$
$$U_\la=u_\la=0\quad\mbox{ on }\partial\Omega,$$
$$\Delta U_\la\in L^1(\Omega)\;\,(\mbox{ since }\,
U_\la\in C^2(\overline\Omega)).$$
In view of Lemma \ref{l} we get $u_\la\leq U_\la$ in $\Omega.$
Thus, a standard bootstrap argument  (see \cite{gt}) implies that 
there exists a
solution $u_{\ep}\in C^{2}(\overline{\Omega})$ of (\ref{epsilon}) such that
$$
u_{\la}\leq u_{\ep}\leq U_{\la}
\quad\mbox{in}\;\;\Omega.
$$
Integrating in (\ref{epsilon}) we obtain
$$\di -\int_{\Omega}\Delta
u_{\ep}dx+K_*\int_{\Omega}g(u_{\ep}+\ep)dx=\la
\int_{\Omega}f(x,u_{\ep})dx.$$
Hence
\begin{equation}\label{treitrei}
\di-\int_{\partial\Omega}\frac{\partial
u_{\ep}}{\partial n}ds+K_*
\int_{\Omega}g(u_{\ep}+\ep)dx\leq M,
\end{equation}
where $M>0$ is a positive constant. Taking into account 
the fact that $\di\frac{\partial
u_{\ep}}{\partial n}\leq 0$ on $\partial\Omega,$
relation
(\ref{treitrei}) yields
$\di K_*\int_{\Omega}g(u_{\ep}+\ep)dx\leq M.$ Since
$u_\ep\leq U_{\la}$ in $\overline\Omega,$ from the last inequality we can conclude that
$\di \int_{\Omega}g(U_{\la}+\ep)dx\leq C,$ for some $C>0.$
Thus, for any compact subset
$\omega\subset\subset\Omega$ we have
$$\di \int_{\omega}g(U_{\la}+\ep)dx\leq C.$$
Letting $\ep\ri 0^+,$ the above relation produces
$\di \int_{\omega}g(U_{\la})dx\leq C.$
Therefore
\begin{equation}\label{treitpatru}
\di \int_{\Omega}g(U_{\la})dx\leq C.
\end{equation}
On the other hand, using \eq{udst} and the hypothese $\int_0^1g(s)ds=+\infty,$ it follows
$$\di \int_{\Omega}g(U_{\la})dx\geq
\int_{\Omega}g(c_2\mbox{dist}\,(x,\partial\Omega))dx=+\infty,$$
which contradicts (\ref{treitpatru}).
Hence, $(1)_{\la}$
has no classical solutions and 
the proof of Theorem \ref{th1} is now complete.
\qed
\medskip

\section{Proof of Theorem \ref{th2}}

Fix $\la>0.$ We first note that $U_\la$ defined in \eq{UU} is a super-solution of $(1)_\la.$ 
We foccuss now on finding a sub-solution $\underline u_\la$ 
such that $\underline u_\la\leq U_\la$ in $\Omega.$

Let $h:[0,\infty)\rightarrow[0,\infty)$ be such
that
\begin{equation}\label{doicinci}
\left\{\begin{tabular}{ll}
$h''(t)=g(h(t)),\quad \mbox{ for all }t>0,$\\
$h>0,\quad \mbox{ in }(0,\infty),$\\
$h(0)=0.$\\
\end{tabular} \right.
\end{equation}
Multiplying by $h'$ in \eq{doicinci} and then integrating over 
$[s,t]$  we have
$$(h')^2(t)-(h')^2(s)=2\int^{h(t)}_{h(s)}g(\tau)d\tau,\quad
\mbox{ for all }\,t>s>0.$$ 
Since $\int^1_0g(\tau)d\tau<\infty,$ from the above equality we deduce
that we can extend $h'$ in origin by taking $h'(0)=0$ and so 
$h\in C^2(0,\infty)\cap C^1[0,\infty).$
Taking into account the fact that $h'$ is increasing and $h''$ is decreasing on $(0,\infty),$  the 
mean value theorem implies that
$$\di\frac{h'(t)}{t}=\frac{h'(t)-h'(0)}{t-0}\geq h''(t),\quad\mbox{ for all }\,t>0.$$
Hence $h'(t)\geq th''(t),$ for all $t>0.$ Integrating in the last inequality we get
\neweq{has}
th'(t)\leq 2h(t),\quad\mbox{ for all }\,t>0.
\endeq

Let $\phi_1$ be the normalized positive
eigenfunction corresponding to the
first eigenvalue $\la_1$ of the problem
$$
 \left\{\begin{tabular}{ll}
$-\Delta u=\la u$ \quad & ${\rm in}\ \Omega,\,$\\
$u=0$ \quad & ${\rm on}\ \partial\Omega\,.$\\
\end{tabular} \right.
$$
It is well known that $\phi_1\in C^2(\overline\Omega).$
Furthermore, by Hopf's maximum principle there exist
$\delta>0$ and $\Omega_0\subset\subset\Omega$ such that
$\di |\nabla \phi_1|\geq\delta$ in $\Omega\setminus\Omega_0.$
Let $M=\max\{1,2K^*\delta^{-2}\},$ where $K^*=\max_{x\in\overline\Omega}K(x).$ 
Since
$$\lim_{{\rm dist}\,(x,\partial\Omega)\ri 0^+}
\Big\{-K^*g(h(\phi_1))+M^a(h')^a(\phi_1)|\nabla\phi_1|^a\Big\}=-\infty,$$
by letting $\Omega_0$ close enough to the boundary of $\Omega$ 
we can assume that
\neweq{omegaz}
-K^*g(h(\phi_1))+M^a(h')^a(\phi_1)|\nabla\phi_1|^a<0\quad\mbox{ in }\;\Omega\setminus\Omega_0.
\endeq
We now are able to show that $\underline{u}_{\la}=Mh(\phi_1)$
is a sub-solution of $(1)_\la$ provided $\la>0$ is sufficiently large.
Using the monotonicity of $g$ and \eq{has} we have
\begin{equation}\label{doinoua}
\begin{tabular}{ll}
$\di -\Delta
\underline{u}_{\la}+K(x)g(\underline{u}_{\la})+|\nabla\underline u_\la|^a=$\\
$\quad\leq -Mg(h(\phi_1))|\nabla\phi_1|^2+\la_1Mh'(\phi_1)\phi_1+
K^*g(Mh(\phi_1))+M^a(h')^a(\phi_1)|\nabla\phi_1|^a$\\
$\quad\leq g(h(\phi_1))(K^*-M|\nabla\phi_1|^2)+\la_1Mh'(\phi_1)\phi_1+
M^a(h')^a(\phi_1)|\nabla\phi_1|^a$\\
$\quad\leq g(h(\phi_1))(K^*-M|\nabla\phi_1|^2)+2\la_1Mh(\phi_1)+
M^a(h')^a(\phi_1)|\nabla\phi_1|^a.$\\
\end{tabular}
\end{equation}
The definition of $M$ and \eq{omegaz} yield
\neweq{bord}
\di -\Delta
\underline{u}_{\la}+K(x)g(\underline{u}_{\la})+|\nabla\underline u_\la|^a
\leq 2\la_1Mh(\phi_1)=2\la_1\underline u_\la\quad\mbox{ in }\;\Omega\setminus\Omega_0.\endeq
Let us choose $\la>0$ such that
\neweq{lambda1}
\di \la\frac{\min_{x\in\overline\Omega_0}f(x,Mh(\|\phi_1\|_\infty))}
{M\|\phi_1\|_\infty}\geq 2\la_1.
\endeq
Then, by virtue of the assumption $(f1)$ and \eq{lambda1} we have
$$\di \la\frac{f(x,\underline u_\la)}{\underline u_\la}\geq 
\la\frac{f(x,Mh(\|\phi_1\|_\infty))}{M\|\phi_1\|_\infty}\geq 2\la_1
\quad\mbox{ in }\Omega\setminus\Omega_0.$$
The last inequality combined with \eq{bord} yield
\neweq{bordfinal}
\di -\Delta
\underline{u}_{\la}+K(x)g(\underline{u}_{\la})+|\nabla\underline u_\la|^a
\leq 2\la_1 \underline u_\la\leq \la f(x,\underline u_\la)\quad\mbox{ in 
}\Omega\setminus\Omega_0.\endeq 
On the other hand, from \eq{doinoua} we obtain
\neweq{interion}
\di -\Delta
\underline{u}_{\la}+K(x)g(\underline{u}_{\la})+|\nabla\underline u_\la|^a
\leq K^*g(h(\phi_1))+2\la_1Mh(\phi_1)+
M^a(h')^a(\phi_1)|\nabla\phi_1|^a\quad\mbox{ in }\Omega_0.
\end{equation}
Since $\phi_1>0$ in $\overline\Omega_0$ and $f$ is positive on 
$\overline \Omega_0\times(0,\infty),$ we may choose $\la>0$ such that
\neweq{lambda2}
\di \la\min_{x\in\overline\Omega_0}
f(x,Mh(\phi_1))\geq
\max_{x\in\overline\Omega_0}\Big\{K^*g(h(\phi_1))+2\la_1Mh(\phi_1)+
M^a(h')^a(\phi_1)|\nabla\phi_1|^a\Big\}.
\endeq
From \eq{interion} and \eq{lambda2} we deduce
\neweq{intfinal}
\di -\Delta
\underline{u}_{\la}+K(x)g(\underline{u}_{\la})+|\nabla\underline u_\la|^a
\leq \la f(x,\underline u_\la)\quad\mbox{ in }\Omega_0.
\endeq 
Now, \eq{bordfinal} together with \eq{intfinal} shows that $\underline u_\la=Mh(\phi_1)$ is a 
sub-solution of $(1)_\la$ provided $\la>0$ satisfy
\eq{lambda1} and \eq{lambda2}. With the same arguments as in the proof of Theorem \ref{th1} and 
using Lemma \ref{l}, one can prove that $\underline u_\la\leq U_\la$ in $\Omega.$ By a standard 
bootstrap argument (see \cite{gt}) we obtain a classical 
solution $u_\la$ such that $\underline u_\la\leq u_\la\leq U_\la$ 
in $\Omega.$   

We have proved that $(1)_\la$ has at least one classical 
solution when $\la>0$ is large. Set
$$ \di A=\{\;\la>0; \mbox{ problem }(1)_\la
\mbox{ has at least one classical solution}\}.$$
From the above arguments we deduce that $A$ is nonempty. Let
$\la^*=\inf A.$   
We claim that if $\la\in A,$ then $(\la,+\infty)\subseteq A.$ 
To this aim, let $\la_1\in A$ and $\la_2>\la_1.$ If
$u_{\la_1}$ is a solution of $(1)_{\la_1},$ then
$u_{\la_1}$ is a sub-solution for $(1)_{\la_2}$ while
$U_{\la_2}$  defined in \eq{UU} for $\la=\la_2$  is a super-solution. 
Moreover, we have 
$$\di \Delta U_{\la_2}+\la_2f(x,U_{\la_2})\leq 0\leq \Delta u_{\la_1}+
\la_2f(x,u_{\la_1})\quad\mbox{ in }\Omega,$$
$$\di U_{\la_2},u_{\la_1}>0\quad\mbox{ in }\Omega,$$
$$U_{\la_2}=u_{\la_1}=0\quad\mbox{ on }\partial\Omega$$
$$\Delta U_{\la_2}\in L^1(\Omega).$$
Again by Lemma \ref{l} we get $u_{\la_1}\leq U_{\la_2}$ in $\Omega.$
Therefore, the problem $(1)_{\la_2}$ has at least one 
classical solution. This proves the claim. Since
$\la\in A$ was arbitrary chosen, we conclude that
$(\la^*,+\infty)\subset A.$ 

To end the proof, it suffices to show that $\la^*>0.$ In that sense, we 
will prove that there exists $\la>0$ small enough such that $(1)_\la$
has no classical solutions. We first remark that
$$\di \lim_{s\ri 0^+}(f(x,s)-K(x)g(s))=-\infty \quad
\mbox{ uniformly for }\, x\in\Omega.$$
Hence, there exists $c>0$ such that
\neweq{n1}
\di f(x,s)-K(x)g(s)<0, \quad
\mbox{ for all }\, (x,s)\in\Omega\times(0,c).
\endeq
On the other hand, the assumption $(f1)$ yields 
\neweq{n2}
\di \frac{f(x,s)-K(x)g(s)}{s}\leq \frac{f(x,s)}{s}\leq
\frac{f(x,c)}{c} \quad
\mbox{ for all }\, (x,s)\in\Omega\times[c,+\infty).
\endeq 
Let $m=\max_{x\in\overline\Omega}\frac{f(x,c)}{c}.$ 
Combinind \eq{n1} with \eq{n2} we find
\neweq{n3}
\di f(x,s)-K(x)g(s)<ms, \quad
\mbox{ for all }\, (x,s)\in\Omega\times(0,+\infty).
\endeq
Set $\la_0=\min\left\{1,\la_1/2m\right\}.$
We show that problem $(1)_{\la_0}$ has no classical solution. Indeed, 
if $u_0$ would be a classical solution of $(1)_{\la_0},$ then, according to \eq{n3}, $u_0$ is a 
sub-solution of
\neweq{n4} 
\left\{\begin{tabular}{ll}
$\di-\Delta u=\frac{\la_1}{2}u$ \quad & ${\rm in}\
\Omega,$\\
$u>0$ \quad & ${\rm in}\ \Omega,$\\
$u=0$ \quad & ${\rm on}\ \partial\Omega.$\\
\end{tabular} \right. 
\endeq
Obvously, $\phi_1$ is a super-solution of \eq{n4} and by Lemma \ref{l} 
we get $u_0\leq \phi_1$ in $\Omega.$ Thus, by standard elliptic arguments, problem \eq{n4} has a 
solution $u\in C^2(\overline\Omega).$
Multiplying by $\phi_1$ in \eq{n4} and then integrating over 
$\Omega$ we have
$$\di -\int_\Omega\phi_1\Delta udx=\frac{\la_1}{2}\int_\Omega 
u\phi_1dx,$$ 
that is,
$$\di -\int_\Omega u\Delta\phi_1dx=\frac{\la_1}{2}\int_\Omega 
u\phi_1dx.$$ 
The above equality yields $\int_\Omega 
u\phi_1dx=0,$ which is clearly a contradiction, since $u$ and $\phi_1$ are positive on $\Omega.$
If follows that problem $(1)_{\la_0}$ has no classical 
solutions which means that $\la^*>0.$
This completes the proof of Theorem \ref{th2}.
\qed

\medskip
{\bf Acknowledgments.} The authors are partially supported by Programme EGIDE-Brancusi between 
University of Craiova and Universit\'e de Picardie Jules Verne in Amiens. M.~Ghergu is also 
partially supported by Grant CNCSIS TD 25/2005.

\end{document}